\documentclass[12pt,a4paper]{amsart}
\usepackage{amssymb}

\input{epsf}
\newtheorem{theorem}{Theorem}[section]
\def\func#1{\mathop{\rm #1}}%

\makeatletter
\def\RIfM@{\relax\protect\ifmmode}
\def\text{\RIfM@\expandafter\text@\else\expandafter\mbox\fi}
\let\nfss@text\text
\def\text@#1{\mathchoice
   {\textdef@\displaystyle\f@size{#1}}%
   {\textdef@\textstyle\tf@size{\firstchoice@false #1}}%
   {\textdef@\textstyle\sf@size{\firstchoice@false #1}}%
   {\textdef@\textstyle \ssf@size{\firstchoice@false #1}}%
   \glb@settings}
\def\textdef@#1#2#3{\hbox{{%
                    \everymath{#1}%
                    \let\f@size#2\selectfont
                    #3}}}
\newif\iffirstchoice@
\firstchoice@true
\makeatother
\def\stackunder#1#2{\mathrel{\mathop{#2}\limits_{#1}}}%

\begin{document}
%% end of leading block
\title{Jumping oscillator}
\author{F. Pugliese, A. M. Vinogradov}

\address{Dipartimento di Matematica e Informatica, Universit\`{a} di Salerno,
Via S. Allende, I-84081 Baronissi (SA),Italy
and
Istituto Nazionale di Fisica Nucleare, sezione di Napoli-Salerno}
\email{pugliese@matna2.dma.unina.it}
\address{The Diffiety Institute, Dipartimento di Matematica e Informatica, Universit\`{a} di Salerno,
Via S. Allende, I-84081 Baronissi (SA),Italy
and
Istituto Nazionale di Fisica Nucleare, sezione di Napoli-Salerno}
\email{vinograd@ponza.dia.unisa.it}
\thanks{%
Partially supported by the italian Ministero dell\'Universit\`{a} e
della Ricerca Scientifica e Tecnologica.}

\begin{abstract}
It is shown that a lagrangian system whose Legendre transformation
degenerates along a hypersurface behaves in a strange manner by jumping from
time to time without any ''visible cause''. In such a jump the system
changes instantaneously its coordinates as well as its momenta. The
mathematical dscription of the phenomenon is based on the theory of impact,
refraction and reflection developed by one of the authors and the
observation that a hamiltonian vector field, understood as a relative one,
can be associated with any lagrangian, degenerated or not. Necessary
elements of the general theory of such systems are reported and a detailed
description of a post-relativistic oscillator showing such a behaviour is
given.
\end{abstract}

\maketitle

\tableofcontents

\section{Introduction.}

In Lagrangian mechanics two, in a sense, extremal situations were widely
studied. One of them, classical, corresponds to systems with nowhere
degenerated Legendre map. In this situation the Legendre map identifies the
Lagrangian dynamics with the corresponding Hamiltonian one. On the contrary,
in the second case the Legendre transformation is supposed to be
everywhere degenerated (and of constant rank). Such a situation is common in
gauge theories and by this reason was studied in several works, starting
with the pioneering paper by Dirac \cite{Dirac}. But in a generic, in the
sense of singularity theory, situation the Legendre map $\mathcal{L}$ is
almost everywhere non-degenerated except a critical hypersurface $S$. Since $%
\mathcal{L}$ maps the tangent bundle $T\left( M\right) $ of the
configuration space $M$ of the system in question to the cotangent bundle $%
T^{*}\left( M\right) $, which both have the same dimension, the most generic
singularities of $\mathcal{L}$ are just folding points (see, for instance,
\cite{Golub-Guill}). In such a case $S$ subdivides $T\left( M\right) $ or
a suitable domain of it into a number of regular regions, inside each of
which the dynamics admits the standard description either in Lagrangian or
in Hamiltonian form. But what happens when a trajectory reaches the singular
hypersurface $S$? In this paper we show that it jumps. More exactly, if its
arrival point is of the fold type , i.e. a generic one, then the moving
point in $M$ at a certain instant disappears suddenly and at the same time
reappears in one or more, generally, distant points to continue its smooth
motion up to an eventual subsequent jump. It is worth to emphasize that
collisions and impacts of any kind in mechanics as well as refractions and
reflections in geometrical optics in the context of the hamiltonian approach
are phenomena of the same mathematical nature. In this case the ''jumping''
quantity is the velocity, but not the space-time coordinates (see
\cite{Vin-Kup},appendix 3).

In this paper, which is the first of a series planned to be written in this
connection, we present the general mathematical background and then
illustrate it with an example of a certain physical flavour, a relativistic
oscillator in the post-Galilean approximation. Our approach is a fusion of
two simple but, it seems, important observations.. The first of them is that
the hamiltonian vector field can be associated naturally with any
lagrangian, independently of the fact whether the corresponding Legendre map
is degenerated or not. This field, however, is a relative one (see below)
with respect to the Legendre map. The second observation is an analogy with
the hamiltonian theory of impact as developed by one of the authors (
\cite{Vin-Kup}, appendix 3).

\noindent It was not our goal here to discuss variational aspects of the
problem. It will be done in a separate purely mathematical paper. We only
remark that variational problems with non- regular lagrangians were already
studied by K. Weierstrass and various authors later on (see, for instance,
the Weierstrass-Erdmann theorem in \cite{Gel-Fom}).

It seems that not too much attention was paid to physical systems the
actions of which are not regular everywhere. We mention here recent works by
G. Vilasi, I. Pavlotsky and their collaborators (\cite{Vil-Pav},
\cite{Las-Pav}), where some particular results concerning the post-Galilean
oscillator and the Darwin two electrons-model were obtained. Our approach
is, however, completely different and leads directly to a complete dynamical
picture.

We hope that the forthcoming study of more realistic models will clarify the
physical content of the proposed mechanism. At the moment its eventual
applications look so attractive that it would be wiser to postpone
speculations.

\medskip

\section{ Hamiltonian Theory of Impact: the Transition Principle.}

Here we describe with the necessary details discontinuous hamiltonian
systems, which will help us to understand the behaviour of dynamical systems
described by singular lagrangians. Namely, the geometry that describes jumps
of phase trajectories for such hamiltonian systems is the same as for
lagrangian ones with singularities of fold type. In other words, the
principle controlling these ''jumps'' in the pure hamiltonian case is quite
transparent and shows what one has to do in the more complicated lagrangian
situation.

Let $(\Phi ,\Omega )$ be the phase space of a dynamical system with $\Omega
=\sum_{i}dp_{i}\wedge dq_{i}$ being a symplectic 2-form on $\Phi $. Suppose
then that $\Phi $ is divided by a hypersurface $\Gamma $ into two closed
domains $\Phi _{+}$,$\Phi _{-}$, having $\Gamma $ as their common boundary,
i.e. $\partial \Phi _{+}=\Gamma =\partial \Phi _{-}$. Suppose also that the
hamiltonian of the system is smooth on $\Phi _{\pm }$. In other words, if H$%
_{\pm }=H|_{\Phi _{\pm }}$, then $H_{\pm }\in C^{\infty }\left( \Phi _{\pm
}\right) $. This, in particular, means that $H_{\pm }|_{\Gamma }\in
C^{\infty }\left( \Gamma \right) $, but it is not supposed that $%
H_{+}|_{\Gamma }$ coincides with $H_{-}|_{\Gamma }$. So $H$, and
consequently the hamiltonian field $X_{H}$ associated to it, is well defined
on $\Phi \backslash \Gamma $ and is bi-valued on $\Gamma $.

As it is well known, the local coordinate expression of $X_{H_{\pm }}$ is
\[
X_{H_{\pm }}=\sum_{i}(\frac{\partial H_{\pm }}{\partial p_{i}}\frac{\partial
}{\partial q_{i}}-\frac{\partial H_{\pm }}{\partial q_{i}}\frac{\partial }{%
\partial p_{i}})\quad \,\,\,\,\,\,\,\,,\,\,\,\,\,\,\,\,\,\,\,\,\,\,\,\,\,\,
\]
The corresponding canonical equations :
\[
\stackrel{\cdot }{q_{i}}=\frac{\partial H_{\pm }}{\partial p_{i}}\quad
,\quad \stackrel{\cdot }{p_{i}}=-\frac{\partial H_{\pm }}{\partial q_{i}}%
\qquad ,
\]
describe the motion of the system \textit{inside }$\Phi _{\pm }$.But when
the phase trajectory arrives at $\Gamma $ it must ''decide'' under control
of which hamiltonian to proceed on. The \textit{transition principle (
\cite{Vin-Kup})} prescribes how this decision should be taken. Below we recall
this ''recipe''.

First, remember that to each point $x\in \Gamma $ the \textit{characteristic
line} $l_{x}\subset T_{x}\left( \Phi \right) $ of $\Gamma $ at $x$ is
associated (\textit{\cite{Vin-Kup}}). This one-dimensional subspace is the
skew-orthogonal complement of the hyperplane $T_{x}(\Gamma )\subset
T_{x}(\Phi )$ with respect to the bilinear, skew-symmetric, non-degenerate
form $\Omega _{x}$. In other words:
\[
l_{x}=\left\{ \xi \in T_{x}(\Phi )\mid \Omega _{x}(\xi ,\eta )=0\quad
\forall \eta \in T_{x}(\Gamma )\right\} \quad .
\]
Since $l_{x}\subset T_{x}(\Gamma )$, for any $x\in \Gamma $ , this way we
get a one-dimensional distribution $x\mapsto l_{x}$ on $\Gamma $ , whose
integral curves are called \textit{characteristics }of $\Gamma $.

{\bf Example}. Let $S$ be a hypersurface in a manifold $M$, $\dim M=n$.
Then , denoting by $\pi :T^{*}\left( M\right) \rightarrow M$ the canonical
projection, we obtain the hypersurface $\Gamma =\pi ^{-1}\left( S\right)
\subset \Phi =T^{*}(M)$. Its characteristics are straight lines contained in
the fibers $T_{q}^{*}\left( M\right) $, $q\in S$. In fact, let $q_{n}=0$ be
the equation of $S$ in a certain local chart. Then this is also the equation
of $\Gamma $ in the corresponding chart $\left( q,p\right) $ on $\Phi $. By
definition a vector $\xi \in T_{x}\left( \Phi \right) ,x\in \Gamma $, is
collinear to the characteristic direction at $x$ iff covector $\Omega
_{x}\left( \xi ,\bullet \right) $vanishes on $T_{x}\left( \Gamma \right) $,
i.e iff $\Omega _{x}\left( \xi ,\bullet \right) =\lambda \,dq_{n}$ for some $%
\lambda \in \mathbf{R}$. It is easy to see that vector $\frac{\partial }{%
\partial p_{n}}|_{x}$ satisfies this condition. Hence characteristics are
the integral curves of the vector field $\frac{\partial }{\partial p_{n}}$,
i.e. straight lines.\medskip

\begin{figure}
\mbox{
\epsfxsize=2.165in
\epsfbox{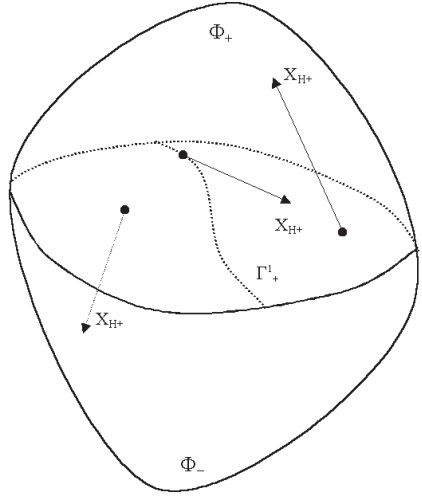}
}
\caption{In- and out-points}
\label{fig1}
\end{figure}

Let $x\in \Gamma $. We say that $x$ is a $+$\textit{-in-point} (resp., a%
\textit{\ }$+$\textit{-out-point}) if $X_{H_{+}}|_{x}$ is directed toward $%
\Phi _{+}$(resp. $\Phi _{-}$). Similarly, we say that $x$ is a $-$\textit{%
-in-point} (resp., a\textit{\ }$-$\textit{-out-point}) if $X_{H_{-}}|_{x}$
is directed toward $\Phi _{-}$(resp. $\Phi _{+}$) (see Fig. \ref{fig1}). In- and
out-points of $H_{\pm }$ are separated by a hypersurface $\Gamma _{\pm
}^{1}\subset \Gamma $ along which $X_{H_{\pm }}$ is tangent to $\Gamma $.

Suppose now that the phase trajectory, starting from a point inside of $\Phi
_{+}$ (resp. $\Phi _{-}$), reaches a point $x\in \Gamma $ at an instant $%
\overline{t}$, and let $E$ be the constant value of $H_{+}$ (resp. $H_{-}$)
along the phase trajectory for $t\leq \overline{t}$ . Denote by $\gamma _{x}$
the characteristic curve of $\Gamma $ passing through $x$, and by $\Sigma
_{E}^{+}$ (resp. $\Sigma _{E}^{-}$) the hypersurface $\{H_{+}=E\}$ (resp., $%
\{H_{-}=E\}$) of $\overline{\Phi }_{+}$ (resp., $\overline{\Phi }_{-}$). A
point $y\in \gamma _{x}\cap $ $\Sigma _{E}^{+}$ (resp. $\gamma _{x}\cap
\Sigma _{E}^{-}$) is called \textit{decisive} for $x$ if it is an $+$%
-in-point (resp. a $-$-in-point). Now we can state the following\smallskip

{\bf Transition Principle}: \textit{When a moving point reaches the
separation hypersurface }$\mathit{\Gamma }$\textit{\ at a point }$x$\textit{%
\ with energy }$E$\textit{, it continues then its motion from all }$x$%
\textit{-decisive points }$y\in \Sigma _{E}^{\pm }$\textit{\ simultaneously
under control of corresponding Hamiltonians }$H_{\pm }$\textit{. The passage
of the phase point from }$x$\textit{\ to the }$y^{\prime }s$\textit{\ is
assumed to be instantaneous.}

\begin{figure}
\mbox{
\epsfxsize=3.943335in
\epsfbox{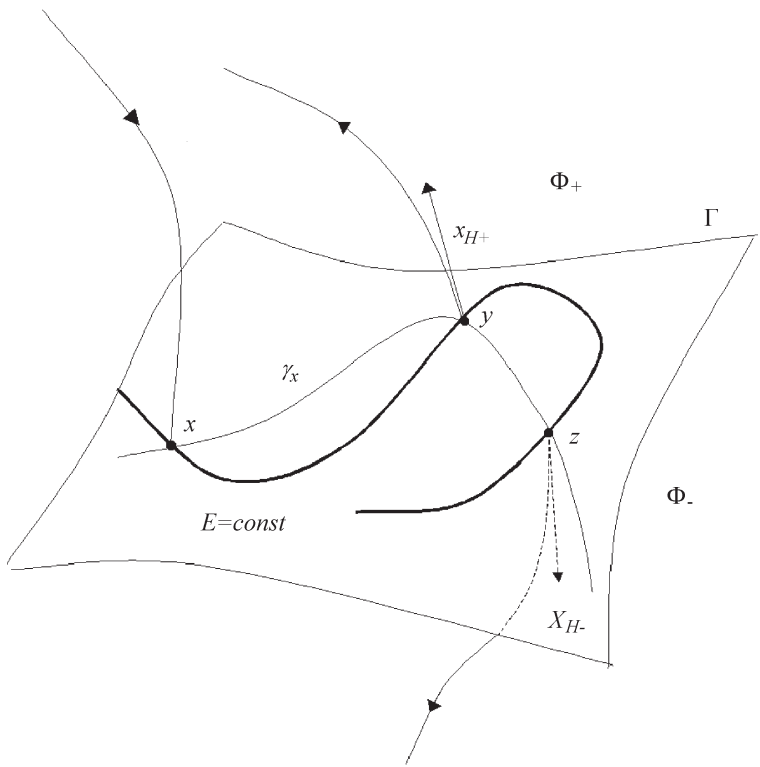}
}
\caption{The Transition Principle ($y$ and $z$ are decisive for $x$)}
\label{fig2}
\end{figure}

The transition principle is illustrated in Fig. \ref{fig2}.

In \cite{Vin-Kup} some examples illustrating this principle in mechanics
and in geometrical optics are given. To these we will now add the following
one.

{\bf Example (laws of reflection and refraction).}

Let $M$ be an inhomogeneous, anisotropic optical medium, referred to a
system of orthogonal coordinates $(q_{1},q_{2},q_{3})$. Denote by $%
V(q)=c/n\left( q\right) $ the light velocity at point $q$. Here $c$ is the
light velocity in the vacuum and $n\left( q\right) $ stands for the
refraction index. Then (\cite{LeviCiv}) the propagation of the light
rays is described by the canonical system:
\[
\stackrel{\cdot }{q_{i}}=H_{p_{i}}\quad ,\quad \stackrel{\cdot }{p_{i}}%
=-H_{q_{i}}\qquad ,\quad i=1,2,3
\]

with
\begin{equation}
H(q,p)=V(q)\left\| p\right\| \,\,\,\,\,\,\,\,\,\,,  \label{ham. optics}
\end{equation}
where $\left\| p\right\| =\sqrt{\stackunder{i}{\sum }p_{i}^{2}}$.

Let now $M_{+},M_{-}$ be two isotropic optical media with refraction indexes
$n_{+}(q)$ and $n_{-}(q)$ respectively, separated by a surface $S$. Suppose
that a ray starting from a point in $M_{+}$ reaches $S$ at a point $%
\overline{q}\equiv \left( \overline{q}_{1},\overline{q}_{2},\overline{q}%
_{3}\right) $ with an impulse $\overline{p}\equiv \left( \overline{p}_{1},%
\overline{p}_{2},\overline{p}_{3}\right) $. The corresponding velocity is $%
\overline{v}=\left( \overline{v}_{1},\overline{v}_{2},\overline{v}%
_{3}\right) $, with $\overline{v}_{i}=\stackrel{.}{q}_{i}=H_{p_{i}}\left(
\overline{q},\overline{p}\right) =V\left( \overline{q}\right) \overline{p}%
_{i}/\left\| \overline{p}\right\| $.

Pose: $M=\overline{M}_{+}\cup \overline{M}_{-}\,,\Phi =T^{*}\left( M\right)
\,,\Phi _{\pm }=T^{*}\left( M_{\pm }\right) $ and $\Gamma =\stackunder{q\in S%
}{\cup }T_{q}^{*}\left( M\right) $. So, the hypersurface $\Gamma $ separates
$\Phi _{+}$ from $\Phi _{-}$ and the phase trajectory of the ray reaches $%
\Gamma $ at a point $\overline{x}\equiv \left( \overline{q},\overline{p}%
\right) $. By the previous example the characteristic curve $\gamma _{%
\overline{x}}$ of $S$ passing through $x$ is a straight line contained in
the fiber $T_{\overline{q}}^{*}\left( M\right) $. Let us choose coordinates $%
\left( q_{1},q_{2},q_{3}\right) $ in such a way that $S$ is tangent at $%
\overline{q}$ to the hyperplane $\left\{ q_{3}=\overline{q}_{3}\right\} $
and the $q_{3}$-axis is directed toward $M_{-}$. Then the characteristic
direction at $\overline{x}$ is parallel to $\frac{\partial }{\partial p_{3}}$
and $\gamma _{\overline{x}}$ is described by equations
\begin{eqnarray}
q_{i}(t) &=&\overline{q}_{i}  \nonumber \\
p_{i}(t) &=&\overline{p}_{i}+\delta _{i3}t\quad ,\quad t\in \mathbf{R}%
,~~i=1,2,3  \label{reflect. char.}
\end{eqnarray}

If $E=H_{+}(\overline{q},\overline{p})$, then, in view of (\ref
{reflect.
char.}) and (\ref{ham. optics}), the intersection $\gamma _{%
\overline{x}}\cap \Sigma _{E}^{+}$ is composed of two points $\overline{x}$
and $x^{*}$ corresponding to values $\overline{t}=0$ and $t^{*}=-2\overline{p%
}_{3}$, respectively. The point $x^{*}\equiv \left( \overline{q}%
,p^{*}\right) $ with $p^{*}=\left( \overline{p}_{1},\overline{p}_{2},-%
\overline{p}_{3}\right) $ is decisive for $\overline{x}$. In fact, by the
choice of coordinates:
\[
\left( X_{H_{+}}\right) _{x^{*}}\left( q_{3}\right) =\frac{\partial H_{+}}{%
\partial p_{3}}\left( x^{*}\right) =-V_{+}(\overline{q})\frac{\overline{p}%
_{3}}{\left\| \overline{p}\right\| }=-\left( X_{H_{+}}\right) _{\overline{x}%
}\left( q_{3}\right) <0\qquad ,
\]
Therefore, the transition principle tells that the reflected ray does always
exist, and starts from the same point $\overline{q}\in S$ in the direction $%
v^{*}=V\left( \overline{q}\right) p^{*}/\left\| p^{*}\right\| $
corresponding to $p^{*}$. Further, $\overline{v}$, $v^{*}$ are coplanar with
the normal to $S$ at $\overline{q}$ and form with it equal angles $\phi
,\psi _{+}$ ,respectively. In fact
\[
\cos \phi =\frac{\left| v_{3}^{*}\right| }{\left\| v^{*}\right\| }=\frac{%
\left| p_{3}^{*}\right| }{\left\| p^{*}\right\| }=\frac{\left| \overline{p}%
_{3}\right| }{\left\| \overline{p}\right\| }=\cos \psi _{+\qquad ,}
\]

On the other hand the intersection $\gamma _{\overline{x}}\cap \Sigma
_{E}^{-}$ corresponds to the values of $t$ satisfying equation
\begin{equation}
t^{2}+2\overline{p}_{3}t+\left( 1-\overline{n}^{2}\right) \left\| \overline{p%
}^{2}\right\| =0\qquad ,  \label{refract.equation}
\end{equation}
with $\overline{n}=n_{-}(\overline{q})/n_{+}(\overline{q})$. Equation $%
\left( \ref{refract.equation}\right) $ admits real solutions iff $\left(
\frac{\overline{p}_{3}}{\left\| \overline{p}\right\| }\right) ^{2}=\left(
\frac{\overline{v}_{3}}{V_{+}}\right) ^{2}\geq 1-\overline{n}^{2}$, i.e. iff
\[
\sin \phi =\sqrt{1-\frac{\overline{v}_{3}^{2}}{V_{+}\left( \overline{q}%
\right) ^{2}}}\leq \overline{n}\qquad .
\]
Therefore, for values of $\phi $ greater than $\arcsin \overline{n}$ there
is no refraction (''total reflection''). If instead $\phi \leq \arcsin
\overline{n}$ (or if $\overline{n}>1$), then $\left( \ref{refract.equation}%
\right) $ has two real solutions, $t=-\overline{p}_{3}\pm \widetilde{p}_{3}$%
, with $\widetilde{p}_{3}=\left\| \overline{p}\right\| \sqrt{\left( \frac{%
\overline{p}_{3}}{\left\| \overline{p}\right\| }+\overline{n}^{2}-1\right) }$%
, to which correspond the two intersections $(\overline{q},\overline{p}_{1},%
\overline{p}_{2},\widetilde{p}_{3})$, $(\overline{q},\overline{p}_{1},%
\overline{p}_{2},-\widetilde{p}_{3})$.Of these, only the first one is
decisive for $\overline{x}$ (one can see it, as before, by considering the
orientation of $X_{H_{-}}$). Therefore, for $0\leq \phi \leq \arcsin
\overline{n}$ the refracted ray does exist, and its (initial) direction $%
\widetilde{p}=\left( \overline{p}_{1},\overline{p}_{2},\widetilde{p}%
_{3}\right) $ is coplanar with the incident ray and the normal to $S$ at $%
\overline{q}$, and forms with this an angle $\psi _{-}$ such that (\textit{%
Snellius' law}):
\[
\frac{\sin \phi }{\sin \psi _{-}}=\frac{\sqrt{1-\frac{\overline{p}_{3}^{2}}{%
\left\| \overline{p}\right\| ^{2}}}}{\sqrt{1-\frac{\widetilde{p}_{3}^{2}}{%
\left\| \widetilde{p}\right\| ^{2}}}}=\overline{n}
\]

\section{Lagrangians with Singular Hypersurfaces.}

In this section the Hamiltonian theory of impact and refraction, as
described in the previous section, is extended to the singular lagrangians.
This will allow us, as it was already mentioned in the introduction, to
describe discontinuities of the motion that occur when the phase point of
the system reaches the singular surface. We start with recalling the concept
of relative vector field.

\subsection{Relative vector fields}

Let $M$ and $N$ be differentiable manifolds, related one another by smooth
mapping $F:M\rightarrow N$. An $\mathbf{R}$-linear operator $X:C^{\infty
}(N)\rightarrow C^{\infty }(M)$ is called a \textit{relative vector field}
on $N$ along $F$ if it satisfies the Leibniz rule:
\[
X(fg)=F^{*}(f)X(g)+F^{*}(g)X(f)\hspace{0pt}\qquad f,g\in C^{\infty }(N)
\]
If $f\in C^{\infty }\left( M\right) $ and $X$ is a relative vector field,
then $fX$ is also a such one. Therefore relative vector fields on $N$ along $%
F$ form a $C^{\infty }\left( M\right) $-module. Denote it by $\mathcal{D}%
\left( N,M;F\right) $. An ''absolute'' vector field $X$ on $N$ can be
considered as a relative one along the identity map $id_{N\text{ }}$.

Geometrically, a relative vector field $X\in \mathcal{D}\left( N,M;F\right) $
can be thought as a map which associates to any point $a\in M$ a vector $%
X_{a}$ tangent to $N$ at the point $F\left( a\right) $. Namely
\[
X_{a}(f):=\left[ X(f)\right] (a)\qquad ,\qquad f\in C^{\infty }(N)
\]

Conversely, to each differentiable map $\sigma :M\rightarrow T(N)$ such that
$p_{N}\circ \sigma =F$ ,where $p_{N}:T(N)\rightarrow N$ is the canonical
projection, one can associate the relative vector field $X$ defined by:
\[
\left[ X(f)\right] (a):=\sigma (a)(f)\quad ,\qquad a\in M,f\in C^{\infty
}(N)
\]

If $x=(x_{1},.....,x_{m})$, $y=(y_{1},......,y_{n})$ are local coordinates
on $M$ and $N$ respectively, the corresponding local expression of $X\in
\mathcal{D}(N,M;F)$ is:
\[
X=\sum_{i=1}^{n}X^{i}(x)(F^{*}\circ \frac{\partial }{\partial y_{i}}),
\]

where functions $X^{i}(x)$ are the components of $X_{x}$ with respect to the
basis $\frac{\partial }{\partial y_{1}}\mid _{F(x)},.......,\frac{\partial }{%
\partial y_{n}}|_{F(x)}$ of $T_{F(x)}(N)$. Below we will use the simplified
notation
\begin{equation}
X=\sum_{i=1}^{n}X^{i}(x)\frac{\partial }{\partial y_{i}}\;\;\;,
\label{relat vecfield}
\end{equation}

having in mind that the vector field (\ref{relat vecfield}) acts on a
function $\phi \left( y\right) \in C^{\infty }\left( N\right) $ as follows:
\[
X\left( \phi \right) \left( x_{1},....,x_{m}\right) =\sum_{i}X_{i}\left(
x\right) \frac{\partial \phi }{\partial y_{i}}\left( y_{1}\left( x\right)
,.....,y_{n}\left( x\right) \right) \,\,\,\,\,\,,
\]

with $y_{i}=y_{i}\left( x\right) $ being the coordinate expression of $F$.

For more details concerning relative vector fields, see \cite{Mod-Vin}.

\subsection[Singular Lagrangians]{Singular Lagrangians and Relative Hamiltonian Theory of Impact
and Refraction}

Let $M$, $\dim M=n$, be the configuration space of a dynamical system
described by a lagrangian $L\in C^{\infty }(T(M))$, and let $\mathcal{L}%
:T(M)\rightarrow T^{*}(M)$ be the corresponding Legendre mapping. Recall
that, in a fixed local chart $(q_{1},....,q_{n})$ on $M$, $\mathcal{L}$ is
represented by equations:
\begin{eqnarray}
q_{i} &=&q_{i}\quad ,\quad \quad i=1,....n  \label{legendre} \\
p_{i} &=&L_{v_{i}}(q,v)\quad ,\quad \quad i=1,....n\qquad ,  \nonumber
\end{eqnarray}

\noindent where $(q,v)$,$(q,p)$ are special coordinates on $T(M)$ and $%
T^{*}(M),$ respectively, associated with $(q_{1},.......,q_{n})$.

Let $S$ be the singular points locus of the Legendre map :
\[
S:=\left\{ x\in T(M)\,|~\mathrm{rk}~d_{x}\mathcal{L}<n\right\}
\]

\noindent This hypersurface can be described as follows. The Jacobian matrix
of $\mathcal{L}$ at point $(q,v)$ is:
\begin{equation}
d_{(q,v)}\mathcal{L\equiv }\left\|
\begin{array}{ll}
\mathbf{1} & \mathbf{0} \\
L_{vq} & L_{vv}
\end{array}
\right\|  \label{differential}
\end{equation}

\noindent where $L_{vq}=\left\| L_{v_{i}q_{j}}\right\| $ and L$_{vv}=\left\|
L_{v_{i}v_{j}}\right\| $ is the hessian matrix . Since the determinant of (%
\ref{differential}) is $\mathcal{H}:=\det L_{vv}$, the equation of $S$ is:
\[
\mathcal{H}(q,v)=0
\]
\noindent We will assume that $S$ is a regular hypersurface of $T(M).$This
assumption implies that for any point $(q,v)\in S$, $d_{(q,v)}\mathcal{H}$
does not vanish.

According to the standard procedure, the motion of the system outside $S$
can be described by Euler-Lagrange equations:
\begin{equation}
\left\{
\begin{array}{c}
\;\;\;\;\;\;\;\;\;\;\,\,\dot{q}_{i}=v_{i} \\
\frac{d}{dt}\left( L_{v_{i}}\right) -L_{q_{i}}=0\;\;,
\end{array}
\right.  \label{Lagr. eq.}
\end{equation}

$i=1,.....,n$. Equations (\ref{Lagr. eq.}) can be rewritten in the normal
form:
\begin{equation}
\stackrel{\cdot }{v_{i}}=f_{i}(q,v)  \label{Lagr. expl.}
\end{equation}

\noindent in a neighbourhood of any point $(q,v)\in T(M)\diagdown S$. On the
other hand, this is no longer possible if $(q,v)\in S$. Namely,
accelerations are undetermined on $S$ and velocities may have
discontinuities.

Recall (see, for instance, \cite{Golub-Guill}) that $x\in T\left( M\right)
$ is a \textit{fold point} of $\mathcal{L}$ if

\begin{eqnarray}
\text{\ \ \ \ \ \ \ \ \ }Ker\,d_{x}\mathcal{L}\cap T_{x}\left( S\right)
&=&\left\{ 0\right\} \;\;\forall x\in S  \label{transversality} \\
\text{\thinspace \ \ \ \ \ \ \ \ \ \ \ \ \ \ \ \ \ \ \ \ \ \ \ \ \ }d_{x}%
\mathcal{H} &\neq &0\;\;\forall x\in S  \label{simplezero}
\end{eqnarray}
Below we restrict our analysis only to fold points, which for a generic $%
\mathcal{L}$ form an open everywhere dense subset in $S$.

In the next section we will see that the lagrangian of the post-galilean
oscillator exhibits also some irregular singularities, which can be studied
as well.

For our purposes it is important to note that the hamiltonian vector field
associated with a lagrangian exists even if the corresponding Legendre map
possesses some singularities. But in such a case this hamiltonian field
\textit{becomes a relative one}. Namely, define the relative vector field $%
X_{H}^{rel}\in D\left( T^{*}\left( M\right) ,T\left( M\right) ;\mathcal{L}%
\right) $ by:
\begin{equation}
X_{H}^{rel}\stackrel{def}{=}\sum_{i}v_{i}\frac{\partial }{\partial q_{i}}%
+\sum_{i}L_{q_{i}}\left( q,v\right) \frac{\partial }{\partial p_{i}}
\label{rel. ham. field}
\end{equation}

If $\mathcal{L}$ is (locally) regular\noindent , i.e. a (local)
diffeomorphism, then
\[
X_{H}=\left( \mathcal{L}^{-1}\right) ^{*}\circ X_{H}^{rel}
\]
is the standard hamiltonian field associated with $L$.This is easily seen
from the standard coordinate expression of a hamiltonian field:
\[
X_{H}=\sum_{i}(H_{p_{i}}\frac{\partial }{\partial q_{i}}-H_{q_{i}}\frac{%
\partial }{\partial p_{i}})
\]

\noindent In fact, in this case
\begin{equation}
H(q,p)=E(q,v(q,p))\quad ,  \label{H=E}
\end{equation}

\noindent where $v=v(q,p)$ is the inverse of the Legendre map and
\[
E\left( q,v\right) \stackrel{def}{=}\sum_{i}v_{i}L_{v_{i}}\left( q,v\right)
-L\left( q,v\right)
\]
is the energy. Keeping in mind (\ref{legendre}) one finds
\begin{eqnarray*}
H_{p_{j}} &=&\sum_{i}(\delta _{ij}v_{i}+p_{i}\frac{\partial v_{i}}{\partial
p_{j}})-\sum_{i}L_{v_{i}}\frac{\partial v_{i}}{\partial p_{j}}=v_{j} \\
H_{q_{j}} &=&\sum_{i}p_{i}\frac{\partial v_{i}}{\partial q_{j}}%
-L_{q_{j}}-\sum_{i}L_{v_{i}}\frac{\partial v_{i}}{\partial q_{j}}=-L_{q_{j}}
\end{eqnarray*}

\noindent Hence, the vector field $Z_{L}$ corresponding to Euler-Lagrange
equations for $L$ is well-defined on $T\left( M\right) \backslash S$ and,
locally, $$
Z_{L}=\sum_{i}\left( v_{i}\frac{\partial }{\partial q_{i}}%
+f_{i}\left( q,v\right) \frac{\partial }{\partial v_{i}}\right),
$$ with $%
f_{i}$'s given by (\ref{Lagr. expl.}). If $U\subset T\left( M\right)
\backslash S$ is such that $\mathcal{L}_{U}\stackrel{def}{=}\mathcal{L}%
|_{U}:U\rightarrow \mathcal{L}\left( U\right) $ is a diffeomorphism, then
the image $X_{U}\in \mathcal{D}\left( \mathcal{L}\left( U\right) \right) $
of $Z_{L}|_{U}$ with respect to $\mathcal{L}_{U}$ is well defined:
\[
X_{U}=\left( \mathcal{L}_{U}^{-1}\right) ^{*}\circ Z_{L}\circ \mathcal{L}%
_{U}^{*}
\]
and, moreover, is a hamiltonian vector field (in $\mathcal{L}\left( U\right)
$) corresponding to the Hamilton function $H_{U}=E\circ \mathcal{L}%
_{U}^{-1}\in C^{\infty }\left( \mathcal{L}\left( U\right) \right) $. If $%
U^{\prime }$ is another regular domain for $\mathcal{L}$ such that $\mathcal{%
L}\left( U^{\prime }\right) =\mathcal{L}\left( U\right) $, then, generally, $%
H_{U^{\prime }}\neq H_{U}$ or, equivalently, $X_{U}\neq X_{U^{\prime }}$ (in
$\mathcal{L}\left( U\right) $). In other words, in $\mathcal{L}\left(
U\right) $ a multi-valued hamiltonian field is well-defined. So, on the
whole, a multi-valued hamiltonian field is defined and its various branches
are matching one another along $\mathcal{L}\left( S\right) $.

Now we pass to describe the analogues of in- and out-points (see section 2)
in the considered context. The following elementary facts from singularity
theory (see, for instance, (\cite{Golub-Guill})) are needed for this
purpose. First, recall the notion of a submersion with folds. Let $P,Q$ be
two manifolds, with $\dim P=n\geq \dim Q=m$. The \textit{first-order jet} $%
\left[ F\right] _{x}^{1}$ of a map $F:P\rightarrow Q$ at a point $x\in P$
may be considered as a triple $\left( x,y,p\right) $, with $y=F\left(
x\right) $ and $p:T_{x}\left( P\right) \rightarrow T_{y}\left( Q\right) $
being the differential of $F$ at $x$. The manifold of all the first-order
jets of maps from $P$ to $Q$ is denoted by $J^{1}\left( P,Q\right) $. The
submanifold $S_{1}\subset J^{1}\left( P,Q\right) $ is composed of all
triples $\left( x,y,p\right) $ such that $\mathrm{rk}\,p=n-1$. The map $%
j^{1}F:P\rightarrow J^{1}\left( P,Q\right) $ sends a point $x\in P$ to $%
\left[ F\right] _{x}^{1}$. $F$ is called a \textit{submersion with folds} if:

\begin{enumerate}
\item[1)]  $j^1F$ is transversal to $S_1;$

\item[2)]  $Ker\ d_xF$ is transversal to $S_1(F)=(j^1F)^{-1}(S_1)$, for
every $x\in S_1(F)$.
\end{enumerate}

\noindent It is easy to see that in the case $n=m$ condition 1) is
equivalent to the fact that the jacobian of $F$ has only simple zeroes along
$S_{1}(F)$.

The local structure of a submersion with folds is described by the following

\begin{theorem}
\label{submersion} Let $F:P\rightarrow Q$ be a submersion with folds and let
$\stackrel{\_}{x}\in S_1(F)$. Then there exist coordinates $(x_1,....,x_n)$
on $P$ and $(y_{1,....,}y_m)$ on $Q$, centered at $\stackrel{\_}{x}$ and $F(%
\stackrel{\_}{x})$ respectively, in terms of which $F$ takes the form:
\begin{eqnarray*}
y_1=x_1 \\
\vdots  \\
y_{m-1}=x_{m-1} \\
y_m=x_m^2\pm x_{m+1}^2\pm .....\pm x_n^2
\end{eqnarray*}
\smallskip
\end{theorem}

$\blacktriangleleft $ see \cite{Golub-Guill} $\blacktriangleright $.

\noindent In the case we are interested in, $P=T(M),Q=T^{*}(M),F=\mathcal{L}%
,\dim P=\dim Q=2n$. Assumptions (\ref{transversality}), (\ref{simplezero})
guarantee $\mathcal{L}$ to be a submersion with folds. According to theorem
\ref{submersion}, its normal form is:
\begin{equation}
\begin{array}{l}
y_1=x_1 \\
.................. \\
y_{2n-1}=x_{2n-1} \\
y_{2n}=x_{2n}^2
\end{array}
\label{normal form}
\end{equation}

\noindent It results easily from the above description and from (\ref
{normal
form}) that the range of $\mathcal{L}$ locally belongs to the
half-space $y_{2n}\geq 0$ . This allows us to extend the definition of in-
and out-points and the corresponding transition principle to the lagrangian
case. Namely, a point $x\in S$ is called an \textit{in-point} if $%
X_{H}^{rel}|_{x}$ is directed toward the range of $\mathcal{L}$,$%
X_{H}^{rel}|_{x}\left( y_{2n}\right) >0$ , while it is called an \textit{%
out-point} if $X_{H}^{rel}|_{x}$ is directed outside of it, $%
X_{H}^{rel}|_{x}\left( y_{2n}\right) <0$.

\noindent Let us note that the pullback $\mathcal{L}^{*}\left( \Omega
\right) $ of the canonical symplectic form on $T^{*}\left( M\right) $ along
the Legendre map is degenerated on $S$ and is of rank $2n-2$ at any fold
point. The restriction $\Omega _{S}\stackrel{def}{=}\mathcal{L}^{*}\left(
\Omega \right) |_{S}$ continues to be of rank $2n-2$ due to (\ref
{transversality}). This means that the kernel $l_{x}$ of $\Omega _{S}$ at a
fold point $x\in S$ is one-dimensional. This way one gets a one-dimensional
distribution on $S_{fold}$. Characteristic curves are integral curves of it.
Denote by $\gamma _{x}$ the characteristic curve passing through $x\in S$.

\noindent Now it is clear how to extend the transition principle to the
lagrangian case. Namely, calling decisive for $x\in S$ any in-point $y\in
\gamma _{x}$ belonging to the same level $\Sigma _{E}$ of energy of $x$, the
principle can be stated as follows.

{\bf Transition Principle (lagrangian case)}. \textit{When a phase point
moving along }$Z_{L}$\textit{\ reaches at an instant a point }$x\in S$%
\textit{, it then continues its motion along all trajectories of }$Z_{L}$%
\textit{\ issuing from points decisive for }$x$\textit{. Moreover, the
passage from }$x$\textit{\ to a decisive point is instantaneous.}

{\bf Remark 1}. Note that the transition principle implies that the
energy \textit{does not change under an impact with }$S$.

{\bf Remark 2}. According to the principle there are in general as many
possible phase trajectories after the impact with $S$ at a point $x$ as are
the points decisive for $x$.

\section{Relativistic Oscillator}

In this section the behaviour of a relativistic oscillator whose lagrangian
possesses fold singularities is analysed on the base of Transition
Principle. We have chosen this example to show the above theory in action
mainly because of its relative simplicity: its phase trajectories and
characteristics can be described analytically without difficulties. However,
we will see that even in this simple case the behaviour of phase
trajectories with respect to the singular surface is rather interesting, at
least from a geometrical point of view.

By applying the transition principle to describe the discontinuities of
motion of relativistic oscillator we discover a rather remarkable
phenomenon. Namely, if the energy exceeds a certain level and at the same
time the velocity is not too high, the oscillator starts jumping. In other
words, after a smooth motion it instantaneously changes its position (as
well as velocity). In the classical example of reflection and refraction of
light, or elastic collision of bodies and particles, there is no
discontinuity of position. It is worth noting that the ''jumping'' motions
of the oscillator are in a good consistency with the smooth ones.

Finally we note that the relativistic oscillator possesses also
singularities of non-fold type. These are very degenerated and the phase
portrait in this region is rather curious.

\subsection{Relativistic Oscillator}

Recall that there were proposed various relativistic generalizations of the
standard harmonic oscillator (see for instance \cite{Vil-Pav}). Some of
them possess singularities, while others not. Below we study the
two-dimensional post-galilean oscillator of tensor rank 2 \ (\cite{Vil-Pav}%
), possessing both fold and not fold type singularities:
\begin{equation}
L=L\left( r,x\right) =-mc^{2}\left[ \sqrt{1-x}+\left( \frac{r}{r_{0}}\right)
^{2}\left( 1+\frac{x}{2}\right) \right]  \label{oscillator}
\end{equation}

Here $m,r_{0}$ are the mass and the characteristic length of the oscillator
respectively, linked by the relation $r_{0}=\sqrt{\frac{m}{k}}c$ (with $k$
being the elastic constant); $r$ is the distance between the oscillating
mass and the elastic force centre; $x=\frac{v^{2}}{c^{2}}$ is the square of
oscillator velocity, measured with respect to the light velocity $c$. If we
fix in the plane of motion a system of orthogonal coordinates ($q_{1},q_{2})$
with the origin at the centre of the force, then obviously:
\[
r=\sqrt{q_{1}^{2}+q_{2}^{2}}\qquad \qquad x=\frac{v_{1}^{2}+v_{2}^{2}}{c^{2}}
\]

Note that
\begin{equation}
0\leq x<1\,\,\,\,\,\,,  \label{xlimit}
\end{equation}

due to the fact that $v^{2}<c^{2}$. Below we refer to $M=\mathbf{R}%
^{2}=\{\left( q_{1},q_{2}\right) \}$ as the configuration space. So the
lagrangian (\ref{oscillator}) is defined in the domain $\mathcal{U}\subset
T\left( M\right) =\{\left( q,v\right) \}=\mathbf{R}^{2}\times \mathbf{R}^{2}$%
, defined as:
\[
\mathcal{U}=\left\{ (q,v)\in T(M)\mid q\in M,\left\| v\right\| <c\right\}
\]
However we will often not distinguish between $\mathcal{U}$ and $T(M).$ A
similar convention will be adopted also for the cotangent bundle $T^{*}(M)$).

Below we will systematically use in $\mathcal{U}$ the system of coordinates $%
\left( r,\phi ,x,\theta \right) $ (or equivalently $\left( r,\phi
,x,u\right) $, with $u=\theta -\phi $), where $\phi $ and $\vartheta $ are
the angle between $q_{1}$-axis and $\mathbf{r}\equiv (q_{1},q_{2})$, and the
angle between $\mathbf{r}$ and the velocity vector\textbf{\ }$\mathbf{v}%
\equiv (v_{1},v_{2})$, respectively.

To simplify general considerations concerning lagrangian (\ref{oscillator})
it is convenient to work with a generic lagrangian of the form
\begin{equation}
L=L\left( r,x\right)  \label{L=L(r,x)}
\end{equation}

The energy function
\begin{equation}
E=v_{1}L_{v_{1}}+v_{2}L_{v_{2}}-L  \label{energy}
\end{equation}

takes the following form for lagrangian (\ref{L=L(r,x)}):
\begin{equation}
E\left( r,x\right) =2xL_{x}-L  \label{E(r,x)=}
\end{equation}
which in the case of oscillator (\ref{oscillator}) becomes:\label{aggiungi?
1}
\[
E(r,x)=mc^{2}\left[ \frac{1}{\sqrt{1-x}}+\left( \frac{r}{r_{0}}\right)
^{2}\left( 1-\frac{x}{2}\right) \right]
\]
Lagrangian (\ref{L=L(r,x)}) admits also another integral of motion, namely
the \textit{angular momentum}:
\begin{equation}
I(r,x,u)=\frac{2}{c^{2}}L_{x}(r,x)\left( q_{1}v_{2}-q_{2}v_{1}\right) =\frac{%
2}{c}r\sqrt{x}L_{x}(r,x)\sin u  \label{angular momentum}
\end{equation}

For lagrangian (\ref{oscillator}) it is specified as:
\[
I(r,x,u)=mcr\sqrt{x}\left[ \frac{1}{\sqrt{1-x}}-\left( \frac{r}{r_{0}}%
\right) ^{2}\right] \sin u
\]
The integral $I$ corresponds, via Noether's theorem, to the infinitesimal
symmetry :
\[
X=-q_{2}\frac{\partial }{\partial q_{1}}+q_{1}\frac{\partial }{\partial q_{2}%
}=\frac{\partial }{\partial \phi }\qquad ,
\]

of lagrangian (\ref{L=L(r,x)}).

\subsection{Singular hypersurface of relativistic oscillator}

The Legendre map $\mathcal{L}$ associated with lagrangian (\ref{L=L(r,x)})
is given by:
\begin{eqnarray}
q_{i} &=&q_{i}  \label{Leg. osc.} \\
p_{i} &=&L_{v_{i}}=\frac{2}{c^{2}}L_{x}v_{i}\qquad ,~i=1,2  \nonumber
\end{eqnarray}
The corresponding jacobian matrix (\ref{differential}) in terms of standard
coordinates $(q,v)$ and $\left( q,p\right) $ in $T\left( M\right) $ and $%
T^{*}\left( M\right) ,$ respectively, has the entries:
\[
L_{v_{i}q_{j}}=\frac{2}{c^{2}}\frac{L_{xr}}{r}v_{i}q_{j}\quad ,\quad
L_{v_{i}v_{j}}=\frac{2}{c^{2}}L_{x}\delta _{ij}+\frac{4}{c^{4}}%
L_{xx}v_{i}v_{j}\qquad i,j=1,2
\]
So, the corresponding Hessian is
\[
\mathcal{H}(q,v)=L_{v_{1}v_{1}}L_{v_{2}v_{2}}-L_{v_{1}v_{2}}^{2}=\frac{4}{%
c^{4}}L_{x}\left( L_{x}+2xL_{xx}\right) =\mathcal{H}\left( r,x\right)
\]
It is easy to see that:
\begin{equation}
L_{x}+2xL_{xx}=E_{x}\,\,\,\,\,,  \label{Ex=...}
\end{equation}
and, therefore,
\[
\mathcal{H}\left( r,x\right) =\frac{4}{c^{4}}L_{x}E_{x}
\]

\begin{figure}
\mbox{
\epsfxsize=2.485in
\epsfbox{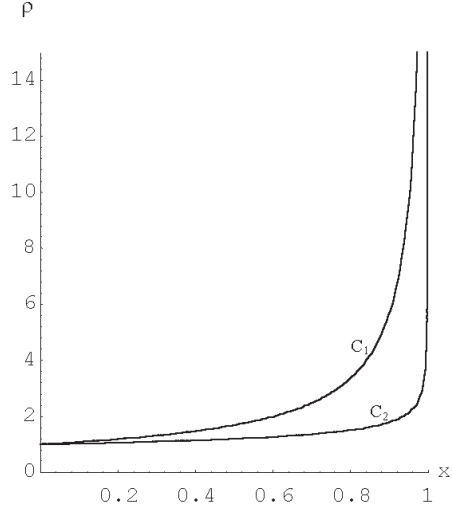}
}
\caption{Singular surface ($r = r/r0$)}
\label{fig3}
\end{figure}

\noindent Hence,
\[
E_{x}L_{x}=0
\]
is the equation of the singular hypersurface $S$. In other words
\[
S=S_{1}\cup S_{2}\qquad ,
\]
with $S_{1}=\{E_{x}=0\},S_{2}=\{L_{x}=0\}$. Each of these hypersurfaces is
fibered in tori $\left( \phi ,\theta \right) $, some of which may reduce to
circles or to a point, depending on $L$. The bases of these ''fibrations''
are curves $C_{1},C_{2}$ in $(x,r)$-plane, given by equations $E_{x}\left(
r,x\right) =0$ and $L_{x}\left( r,x\right) =0,$ respectively. These curves
for the oscillator (\ref{oscillator}) are shown in Fig \ref{fig3}. In this case the
intersection $\gamma =S_{1}\cap S_{2}$ is the circle $\{r=r_{0},x=0\}$
included in the null section $M\subset T(M)$. $\gamma $ is the locus $\func{%
Si}ng~S$ of singular points of $S$. This follows easily from the relations
\begin{eqnarray*}
\mathcal{H}_{q_{i}} &=&\frac{q_{i}}{r}\mathcal{H}_{r}=\frac{q_{i}}{r}\left(
L_{xr}E_{x}+L_{x}E_{xr}\right) \\
\mathcal{H}_{v_{i}} &=&2\frac{v_{i}}{c^{2}}\mathcal{H}_{x}=2\frac{v_{i}}{%
c^{2}}\left( L_{xx}E_{x}+L_{x}E_{xx}\right) \,\,\,\,\,,
\end{eqnarray*}
and from the fact that $S_{1}\backslash \gamma $ and $S_{2}\backslash \gamma
$ are regular.

Now we pass to describe how the kernel of the Legendre map behaves along the
singular hypersurface $S$. It follows from (\ref{differential}) that a
vector $\xi =\sum_{i=1}^{2}\left( a_{i}\frac{\partial }{\partial q_{i}}+b_{i}%
\frac{\partial }{\partial v_{i}}\right) $ belongs to the kernel of $d%
\mathcal{L}$ iff
\begin{equation}
a_{1}=a_{2}=0\qquad ,\quad L_{vv}\left(
\begin{array}{l}
b_{1} \\
b_{2}
\end{array}
\right) =\left(
\begin{array}{l}
0 \\
0
\end{array}
\right)  \label{kernel}
\end{equation}
In full details, the second of conditions (\ref{kernel}) looks as:
\begin{eqnarray}
\left( L_{x}+\frac{2}{c^{2}}L_{xx}v_{1}^{2}\right) b_{1}+\frac{2}{c^{2}}%
L_{xx}v_{1}v_{2}b_{2} &=&0  \label{kernel osc.} \\
\frac{2}{c^{2}}L_{xx}v_{1}v_{2}b_{1}+\left( L_{x}+\frac{2}{c^{2}}%
L_{xx}v_{2}^{2}\right) b_{2} &=&0  \nonumber
\end{eqnarray}
By construction these equations are linearly dependent at any point of $S$.

For a point $\left( q,v\right) \in S_{2}\backslash S_{1}$ there are two
possibilities: $L_{xx}\left( q,v\right) =0$ or $L_{xx}\left( q,v\right) \neq
0$. In the first case the system becomes trivial and $Ker\,d_{\left(
q,v\right) }\mathcal{L}$ is two-dimensional. This never happens for the
oscillator (\ref{oscillator}). In the second case system (\ref{kernel osc.})
reduces to the linear equation:
\[
v_{1}b_{1}+v_{2}b_{2}=0
\]

Except for the points $\left( v_{1}=v_{2}=0\right) $, $Ker\,d_{\left(
q,v\right) }\mathcal{L}$ is one-dimensional. It is generated by the vector $%
v_{2}\frac{\partial }{\partial v_{1}}-v_{1}\frac{\partial }{\partial v_{2}}$%
, which coincides with $-\frac{\partial }{\partial u}$ in coordinates $%
\left( r,\phi ,x,u\right) $. This shows that in both cases $Ker\,d_{\left(
q,v\right) }\mathcal{L}$ is tangent to $S_{2}$. So fold type singularities
do not belong to $S_{2}$ and $\dim \,\mathcal{L}\left( S_{2}\right) <3$. For
the oscillator (\ref{oscillator})$\,\mathcal{L}\left( S_{2}\right) $ is
2-dimensional and coincides with the null section of $T^{*}\left( M\right) $.

Now we go to describe fold points belonging to $S_{1}\backslash S_{2}$.
These form an open domain $S_{1}^{fold}\subset S_{1}$ everywhere dense in $%
S_{1}.$

If $(q,v)\in S_{1}\backslash S_{2}$, it follows from ($\ref{Ex=...})$ that $%
L_{x}=-2xL_{xx}$, hence (\ref{kernel osc.})$_{2}$ becomes
\[
v_{2}b_{1}-v_{1}b_{2}=0
\]
Therefore we have that $Ker~d_{\left( q,v\right) }\mathcal{L}=Span~(v_{1}%
\frac{\partial }{\partial v_{1}}+v_{2}\frac{\partial }{\partial v_{2}}%
)=Span~(\frac{\partial }{\partial x})$ on $S_{1}\backslash S_{2}$.

Thus, singularities of $\mathcal{L}$ along $S_{1}$ are of a substantially
different nature from those along $S_{2}$. Namely, the kernel of $d\mathcal{L%
}$ is \textit{transversal} to $S_{1}\backslash S_{2}$ (due to the fact that $%
E_{xx}\neq 0$ on it) and is \textit{tangent} to $S_{2}$, since $\frac{%
\partial L_{x}}{\partial u}=0$ on it. Hence, $\dim \mathcal{L}\left(
S_{2}\right) <3$ and the transition principle can not be applied to $S_{2}$.
In the next section we shall see that even characteristic directions are
undetermined on $S_{2}$.

\subsection{Characteristic curves on $S$}

To simplify computations we will make use of coordinates $(r,x,\phi
,\vartheta )$ (or, equivalently, $(r,x,u,\phi )$).Let
\[
\rho =\sum_{i}p_{i}\,dq_{i}
\]
be the universal 1-form on $T^{*}(M)$ (\cite{Vin-Kup}). Then
\begin{equation}
\mathcal{L}^{*}(\rho )=L_{v_{1}}dq_{1}+L_{v_{2}}dq_{2}=\frac{2}{c}\sqrt{x}%
L_{x}(r,x)(\cos u\,dr+r\,\sin u\,d\phi )  \label{pullback rho}
\end{equation}
It follows from ($\ref{pullback rho}$) that $\mathcal{L}^{*}(\rho )\mid
_{S_{2}}=0$. Hence
\[
\mathcal{L}^{*}(\Omega )\mid _{S_{2}}=0
\]
Therefore, $\mathcal{L}\left( S_{2}\right) $ is a lagrangian submanifold in $%
T^{*}\left( M\right) $ (with possible singularities) and characteristic
directions on $S_{2}$ are undetermined.

Now we pass to describe characteristics on $S_{1}^{fold}$. Since $grad\,%
\mathcal{H}\neq 0$ on it, the equation $E_{x}(r,x)=0$ of $S_{1}$ can be
solved with respect to one of the variables, say $r$:
\[
r=r_{1}\left( x\right) \,\,\,\,\,\,\,\text{on}\,\,S_{1}^{fold}
\]
In the case of the oscillator (\ref{oscillator}):
\[
r_{1}\left( x\right) =\frac{r_{0}}{(1-x)^{\frac{3}{4}}}
\]
So $(x,u,\phi )$ can be taken as local coordinates on $S_{1}^{fold}$. Then
from (\ref{pullback rho}) we get:
\[
\mathcal{L}^{*}(\rho )\mid _{S_{1}}=\frac{2}{c}\sqrt{x}L_{x}\left(
r_{1}(x),x\right) \left[ r_{1}^{^{\prime }}(x)\cos u\,dx+r_{1}(x)\sin
u\,d\phi \right] \,\,\,\,\,,
\]
so that
\begin{multline}
\mathcal{L}^{*}(\Omega )\mid _{S_{1}}=d\mathcal{L}^{*}(\rho )\mid _{S_{1}}=\\%
\frac{2}{c}\left[ \alpha (x,u)dx\wedge du+\beta (x,u)dx\wedge d\phi +\gamma
(x,u)du\wedge d\phi \right] \,\,\,,  \label{omegaS1}
\end{multline}
with
\begin{eqnarray*}
\alpha (x,u) &=&\sqrt{x}L_{x}(r_{1}(x),x)r_{1}^{^{\prime }}(x)\sin u \\
\beta (x,u) &=&\frac{d}{dx}\left[ r_{1}(x)\sqrt{x}L_{x}\left(
r_{1}(x),x\right) \right] \sin u \\
\gamma (x,u) &=&r_{1}(x)\sqrt{x}L_{x}(r_{1}(x),x)\cos u
\end{eqnarray*}

Characteristic directions on $S_{1}$ are described by a characteristic
vector field $X\in \mathcal{D}\left( S_{1}\right) $, i.e. such that
\begin{equation}
\mathcal{L}^{*}(\Omega )\mid _{S_{1}}(X,\bullet )=0
\label{characteristic condition}
\end{equation}
If
\[
X=a\frac{\partial }{\partial x}+b\frac{\partial }{\partial u}+c\frac{%
\partial }{\partial \phi }\quad ,\qquad a,b,c\in C^{\infty }(S_{1})~,
\]
then (\ref{characteristic condition}) is equivalent, in view of (\ref
{omegaS1}), to
\[
\left(
\begin{array}{lll}
0 & -\alpha & -\beta \\
\alpha & 0 & -\gamma \\
\beta & \gamma & 0
\end{array}
\right) \left(
\begin{array}{l}
a \\
b \\
c
\end{array}
\right) =\left(
\begin{array}{l}
0 \\
0 \\
0
\end{array}
\right)
\]
This system is of rank 2 on $S_{1}^{fold}$, and its fundamental solution is
\[
\left( a,b,c\right) =\left( \gamma ,-\beta ,\alpha \right)
\]

Therefore,
\begin{multline*}
X=\gamma \frac{\partial }{\partial x}-\beta \frac{\partial }{\partial u}%
+\alpha \frac{\partial }{\partial \phi }=\\
q(x)\cos u\frac{\partial }{\partial
x}-q^{^{\prime }}(x)\sin u\frac{\partial }{\partial u}+r_{1}^{^{\prime }}(x)%
\sqrt{x}L_{x}(r_{1}(x),x)\sin u\frac{\partial }{\partial \phi }\quad ,
\end{multline*}

with $q(x)=r_{1}(x)\sqrt{x}L_{x}\left( r_{1}(x),x\right) $ and
characteristic curves are solutions of the system:
\begin{eqnarray}
\stackrel{.}{x} &=&q\left( x\right) \cos u  \nonumber \\
\stackrel{.}{u} &=&-q^{^{\prime }}\left( x\right) \sin u
\label{char. equations} \\
\stackrel{.}{\phi } &=&q\left( x\right) \frac{r_{1}^{^{\prime }}\left(
x\right) }{r_{1}\left( x\right) }\sin u  \nonumber
\end{eqnarray}

Integration of system (\ref{char. equations}) is reduced, obviously, to its
subsystem ($\ref{char. equations})_{1,2}$, whose solutions are to be
described in the rectangle $[0,1[\times \left[ -\pi ,\pi \right] $ due to
the ciclicity of $u$.

It follows from ($\ref{char. equations}$)$_{1,2}$ that
\begin{equation}
\frac{du}{dx}=-\frac{q^{^{\prime }}\left( x\right) }{q\left( x\right) }\tan
u\qquad ,  \label{single char. equat.}
\end{equation}

and, consequently,
\begin{eqnarray*}
\int \frac{du}{\tan u} &=&\ln \left| \sin u\right| =-\int \frac{q^{^{\prime
}}\left( x\right) }{q\left( x\right) }dx=-\int \frac{dq}{q}=-\ln \left|
q\left( x\right) \right| +const. \\
&=&\ln \frac{c}{\left| q\left( x\right) \right| }\qquad ,\quad c>0
\end{eqnarray*}
Hence the general integral of (\ref{single char. equat.}) is
\begin{equation}
\sin u=\frac{a}{q\left( x\right) }\qquad ,\quad a\in \mathbf{R\qquad }
\label{char. curves}
\end{equation}

For the oscillator (\ref{oscillator}) it is specified as:
\[
\sin u=-\frac{2a}{mc^{2}r_{0}}\frac{\left( 1-x\right) ^{\frac{9}{4}}}{x^{%
\frac{3}{2}}}
\]

Let us remark that $q\left( x\right) \sin u=\frac{c}{2}~I\mid _{S_{1}}$. So (%
\ref{char. curves}) shows that $I$ is constant along characteristic curves
of $S_{1}$. Therefore, by the transition principle, angular momentum (as
well as energy) \textit{does not change after the impact with }$S_{1}$.

\begin{figure}
\mbox{
\epsfxsize=2.29in
\epsfbox{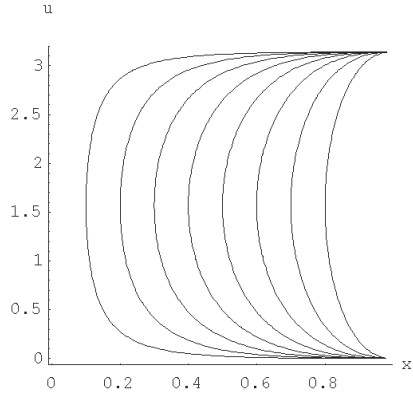}
}
\caption{Characteristic curves on $S1$}
\label{fig4}
\end{figure}

The curves (\ref{char. curves}), denote them by $\gamma _{a}$, for the
oscillator are shown in Fig. \ref{fig4}. Since the variable $u$ is cyclic mod $2\pi $
and $\gamma _{a}$ and $\gamma _{-a}$ are symmetric with respect to the $x$%
-axis we can limit ourselves to dealing with the curves in the rectangle $%
\left( x,u\right) \in [0,1[\times [0,\pi ]$.

\begin{figure}
\mbox{
\epsfxsize=2.818335in
\epsfbox{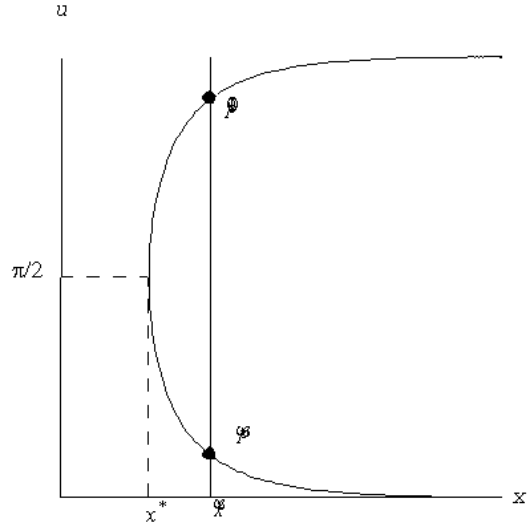}
}
\caption{Intersection between a chracteristic curve and the energy level surface}
\label{fig5}
\end{figure}
Let $\overline{P}\equiv (\bar{x},\bar{u},\bar{\phi})\in S_{1}$ and let $%
\gamma _{\overline{a}}$ be the characteristic passing through $\bar{P}$, $\,%
\bar{E}=E\left( \bar{P}\right) $. The intersection between the energy level
surface $\Sigma _{\bar{E}}$ and $S_{1}$ is the torus $T_{\bar{E}}=\left\{
r=r_{1}\left( \bar{x}\right) ,x=\bar{x}\right\} $. The projection of $T_{%
\overline{E}}$ onto the $\left( x,u\right) $-plane is the line $x=\overline{x%
\text{.}}$ For the oscillator (\ref{oscillator}) this is shown , together
with the projection of $\gamma _{\overline{a}}$, in Fig. \ref{fig5}. Therefore,
assuming $\bar{u}\in [0,\pi ]$, $\gamma _{\bar{a}}$ intersects $T_{\bar{E}}$
at $\bar{P}$ and at $\widetilde{P}\equiv (\bar{x},\pi -\bar{u},\tilde{\phi})$.
In order to determine $\widetilde{\phi }$ notice that
\begin{multline}\label{e:dphi}
\frac{d\phi }{dx}=\frac{r_{1}^{^{\prime }}\left( x\right) }{r_{1}\left(
x\right) }\tan u=\\
\pm \frac{r_{1}^{^{\prime }}\left( x\right) }{r_{1}\left(
x\right) }\frac{\sin u}{\sqrt{1-\sin ^{2}u}}=
\pm \left| \bar{a}\right| \frac{%
r_{1}^{^{\prime }}\left( x\right) }{r_{1}\left( x\right) \sqrt{q^{2}\left(
x\right) -\bar{a}^{2}}}\quad,
\end{multline}
as it results directly from ($\ref{char. equations}$)$_{1,2}$ and (\ref
{char. curves}). In (\ref{e:dphi}) the choice of ''+'' (resp. ''-'')
corresponds to $u\in [0,\frac{\pi }{2}[$ (resp. $u\in ]\frac{\pi }{2},\pi ]$). Hence, the possible position jump is described by the following formula:
\begin{equation}
\widetilde{\phi }-\bar{\phi}=\Delta \phi \left( \bar{x},\bar{u}\right) =\pm
2\left| \bar{a}\right| \int_{x^{*}}^{\bar{x}}\frac{r_{1}^{^{\prime }}\left(
x\right) }{r_{1}\left( x\right) \sqrt{q^{2}\left( x\right) -\bar{a}^{2}}}%
dx\quad \quad ,  \label{deltaphi}
\end{equation}
where $x^{*}$ (see Fig. \ref{fig5}) is the root of the equation
\begin{equation}
\left| q\left( x\right) \right| =\left| \bar{a}\right| \quad \quad \quad ,
\label{equat. for x*}
\end{equation}
and the sign $+$ (resp. $-$) corresponds to $\bar{u}\in \left[ \pi /2,\pi
\right] $ (resp. $\bar{u}\in \left[ 0,\pi /2\right] $). In the case of the
oscillator equation (\ref{equat. for x*}) becomes
\[
x^{2}=\left( \frac{2\left| \overline{a}\right| }{mc^{2}r_{0}}\right)
^{4/3}\left( 1-x\right) ^{3}\quad \quad ,
\]
which has only one root $x^{*}$ in the interval $[0,1[$.
Relation (\ref{e:dphi}) remains valid also for $\bar{u}\in \left[ -\pi ,0\right] $. In
this case the sign + (resp. -) corresponds to $\bar{u}\in $ $\left[ -\pi
,-\pi /2\right] $ (resp. $\bar{u}\in \left[ -\pi /2,0\right] $).

Since on a given characteristic only two points lie, $\overline{P}$ and $%
\widetilde{P}$, belonging to the same energy level, a jump from $\overline{P}
$ to $\widetilde{P}$ or vice-versa may happen only if one of these points is
''in'' while the other is ''out''. This occurs iff the function $%
X_{H}^{rel}\left( g\right) $, $g\left( q,p\right) =0$ being the equation of $%
\mathcal{L}\left( S_{1}\right) $, takes opposite signs at points $\overline{P%
}$, $\widetilde{P}$, and we go to analyse when such is the case.

\begin{figure}
\mbox{
\epsfxsize=4.7666666667in
\epsfbox{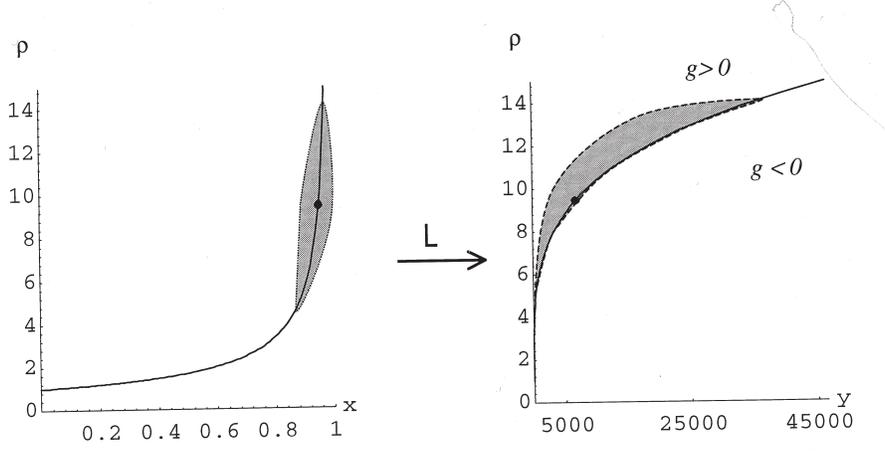}
}
\caption{Legendre map for the oscillator ($* = r/r0$)}
\label{fig6}
\end{figure}

It follows from (\ref{rel. ham. field}) that for lagrangian (\ref{L=L(r,x)})
\begin{equation}
X_{H}^{rel}=v_{1}\frac{\partial }{\partial q_{1}}+v_{2}\frac{\partial }{%
\partial q_{2}}+\frac{L_{r}}{r}(q_{1}\frac{\partial }{\partial p_{1}}+q_{2}%
\frac{\partial }{\partial p_{2}})  \label{ham. field oscill.}
\end{equation}
Below we work with the local chart $\left( r,\phi ,y,\alpha \right) $, $%
y=\left( p_{1}^{2}+p_{2}^{2}\right) /m^{2}c^{2},\,\alpha =\arctan \frac{p_{2}%
}{p_{1}}$, on $T^{*}\left( M\right) $. In terms of these coordinates the
Legendre map (\ref{Leg. osc.}) is given as follows:
\begin{equation}
\quad y=\psi \left( r,x\right) \,,\quad \alpha =\theta \;\;\;\;,
\label{legendre osc.2}
\end{equation}
with $\psi \left( r,x\right) =\frac{4}{m^{2}c^{4}}xL_{x}^{2}\left(
r,x\right) $, and the expression (\ref{ham. field oscill.}) takes the form
\begin{equation}
X_{H}^{rel}=c\sqrt{x}\left( \cos u\frac{\partial }{\partial r}+\frac{\sin u}{%
r}\frac{\partial }{\partial \phi }\right) +\frac{\sqrt{\psi }}{mc}%
L_{r}\left( 2\cos u\frac{\partial }{\partial y}-\frac{\sin u}{\psi }\frac{%
\partial }{\partial \alpha }\right) \quad ,  \label{ham. field oscill.2}
\end{equation}
Due to (\ref{legendre osc.2})$_{1}$ hypersurface $\mathcal{L}\left(
S_{1}\right) $ is given by:
\[
g\left( r,y\right) =0\quad \quad ,
\]
with
\begin{equation}
g\left( r,y\right) =\frac{4}{m^{2}c^{4}}x_{1}\left( r\right) \widetilde{L}%
_{x}^{2}\left( r\right) -y\quad ,  \label{g(r,y)=...}
\end{equation}
$\,$where $x=x_{1}\left( r\right) $ is the function implicitly defined by
equation $E_{x}\left( r,x\right) =0$, and $\widetilde{L}_{x}\left( r\right)
=L_{x}\left( r,x_{1}\left( r\right) \right) $. For the oscillator (\ref
{oscillator}) the Legendre mapping $\mathcal{L}$ is illustrated in Fig \ref{fig6}.
In this case, as it is easy to see from (\ref{legendre osc.2}) and (\ref
{g(r,y)=...}), the image with respect to $\mathcal{L}$ of a sufficiently
small neighbourhood of a point $\left( q,v\right) \in S_{1}$ is contained in
the region $\{g>0\}$.

From (\ref{ham. field oscill.2}) we get:
\[
X_{H}^{rel}\left( g\right) =\cos u\left( c\sqrt{x}g_{r}+\frac{2\sqrt{\psi }}{%
mc}L_{r}g_{y}\right)
\]
But from (\ref{Ex=...}), (\ref{g(r,y)=...}) and $\mathcal{L}\left(
S_{1}\right) $'s equation follows:
\begin{multline*}
g_{r}=\frac{4}{m^{2}c^{4}}\left[ x_{1}^{^{\prime }}\widetilde{L}%
_{x}^{2}+2x_{1}\widetilde{L}_{x}\left( \widetilde{L}_{xr}+x_{1}^{^{\prime }}%
\widetilde{L}_{xx}\right) \right] =\\
\frac{4\widetilde{L}_{x}}{m^{2}c^{4}}%
\left( x_{1}^{^{\prime }}\widetilde{E}_{x}+2x_{1}\widetilde{L}_{xr}\right) =%
\frac{8}{m^{2}c^{4}}x_{1}\widetilde{L}_{x}\widetilde{L}_{xr}
\end{multline*}
where, as before, $\widetilde{f}\left( r\right) $ means $f\left(
r,x_{1}\left( r\right) \right) $. Hence, we obtain:
\begin{equation}
X_{H}^{rel}\left( g\right) =\frac{4}{m^{2}c^{3}}\sqrt{x}\cos u\left( 2x_{1}%
\widetilde{L}_{x}\widetilde{L}_{xr}-L_{x}L_{r}\right)  \label{Xh(g)=...}
\end{equation}
Since $2xL_{xr}-L_{r}=E_{r}$ the restriction of (\ref{Xh(g)=...}) to $S_{1}$
is:
\begin{equation}
\widetilde{X_{H}^{rel}\left( g\right) }=\frac{4}{m^{2}c^{3}}\sqrt{x_{1}}\cos
u\widetilde{L}_{x}\widetilde{E}_{r}  \label{Xh(g)tilde}
\end{equation}

In the case of the oscillator we have $\widetilde{L}_{x}<0,\,\widetilde{E}%
_{r}>0\,$. Therefore (\ref{Xh(g)tilde}) shows that $\widetilde{%
X_{H}^{rel}\left( g\right) }$ is:

1) positive in the region $S_{1}^{+}$ corresponding to $u\in [\pi /2,\pi
]\cup [-\pi ,-\pi /2]$;

2) negative in the region $S_{1}^{-}$ corresponding to $u\in [-\pi /2,\pi
/2] $;

3) null on the bidimensional surface $W=S_{1}\cap \{u=\pm \pi /2\}$.

Hence, keeping in mind what was said above about the range of $\mathcal{L}$
and applying the transition principle, we see that:

1) If $\bar{P}\in S_{1}^{+}$, then $\overline{P}$ is an in-point, while $%
\widetilde{P}$ is an out-point. Therefore, if the phase point, starting from
outside $S_{1}$, reaches it at $\widetilde{P}$, its trajectory can be
prolongated starting from $\bar{P}$ (jump from $\widetilde{P}$ to $\overline{%
P}$).

2)If $\bar{P}\in S_{1}^{-}$, then $\overline{P}$ is an out-point, while $%
\widetilde{P}$ is an in-point. Therefore, if the phase point, starting from
outside $S_{1}$, reaches it at $\overline{P}$, its trajectory can be
prolongated starting from $\widetilde{P}$ (jump from $\overline{P}$ to$%
\widetilde{P}$ ).

3)If $\bar{P}\in W$, i.e. if $\bar{u}=\pm \pi /2$, then $\bar{P}\equiv
\widetilde{P}$. In this case the jump becomes infinitesimal and its
direction is indicated by the hamiltonian vector field. In fact $\mathcal{L}%
\left( W\right) $ is described in terms of coordinates $\left( r,\phi
,y,\alpha \right) $ by equations
\begin{eqnarray*}
g\left( r,y\right) &=&0 \\
\cos u &=&0
\end{eqnarray*}
But it follows from (\ref{ham. field oscill.2}) that on $W$%
\[
X_{H}^{rel}\left( \cos u\right) =-\frac{L_{r}}{mc\sqrt{\psi }}\sin u=\pm
\frac{L_{r}}{mc\sqrt{\psi }}\;\;\;,
\]
which is different from zero on $S_{1}.$ More precisely, it is positive for $%
u=\pi /2$ and negative for $u=-\pi /2$; in both cases $X_{H}^{rel}$ is
directed toward the region S$_{1}^{-}$ in which singular trajectories end.

Finally note that, for the oscillator (\ref{oscillator}),
\[
X_{H}^{rel}\left( p_{i}\right) =\frac{q_{i}}{r}L_{r}\neq 0\;\;\;\text{on}%
\,\,S_{2}\backslash S_{1}
\]
Since $\mathcal{L}\left( S_{2}\right) $ coincides with the null section of $%
T^{*}\left( M\right) $, this shows that $X_{H}$ is transversal to $\mathcal{L%
}\left( S_{2}\right) $.

\subsection{Phase trajectories of the oscillator.}

In this concluding subsection we will study phase trajectories of the
oscillator outside $S$. It will be shown that their behaviour depends
strongly on their position with respect to the singular surface.
Trajectories arriving at $S$ have discontinuities, described by the
transition principle, and together with others form a perfectly
self-consistent dynamical model.

In this subsection coordinates $\left( r,x,u,\phi \right) $ are used. Let,
as before
\[
Z_{L}=\dot{r}\frac{\partial }{\partial r}+\dot{\phi}\frac{\partial }{%
\partial \phi }+\dot{x}\frac{\partial }{\partial x}+\dot{u}\frac{\partial }{%
\partial u}
\]
be the vector field on $T\left( M\right) $ corresponding to the lagrangian $%
L $. Since $I,E$ are first integrals, then
\[
Z_{L}\left( I\right) =Z_{L}\left( E\right) =0\quad ,
\]
or, equivalently
\[
\left(
\begin{array}{lll}
I_{r} & I_{x} & I_{u} \\
E_{r} & E_{x} & 0
\end{array}
\right) \left(
\begin{array}{l}
\dot{r} \\
\dot{x} \\
\dot{u}
\end{array}
\right) =\left(
\begin{array}{l}
0 \\
0 \\
0
\end{array}
\right)
\]

Therefore
\begin{equation}
\dot{r}=-kI_{u}E_{x}\,\,,\quad \dot{x}=kI_{u}E_{r}\,\,,\quad \dot{u}=k\left(
I_{r}E_{x}-I_{x}E_{r}\right)  \label{Zl system}
\end{equation}
On the other hand
\[
\dot{r}=\frac{d}{dt}\left( \sqrt{q_{1}^{2}+q_{2}^{2}}\right) =\frac{%
q_{1}v_{1}+q_{2}v_{2}}{r}=c\sqrt{x}\cos u
\]
So that
\begin{equation}
k=-\frac{c\sqrt{x}\cos u}{E_{x}I_{u}}=-\frac{c\sqrt{x}\sin u}{IE_{x}}=-\frac{%
c^{2}}{2rL_{x}E_{x}}  \label{k=...}
\end{equation}
Notice also the relation:
\begin{equation}
I_{r}E_{x}-I_{x}E_{r}=\frac{E_{x}}{c\sqrt{x}}\sin u\left(
2xL_{x}+rL_{r}\right) \,\,\,\,\,\,\,\,\,\,\,,  \label{temporary1}
\end{equation}
which follows directly from (\ref{E(r,x)=}), (\ref{Ex=...}), (\ref
{angular
momentum}), and the relation
\[
\frac{\partial }{\partial x}\left( \sqrt{x}L_{x}\right) =\frac{E_{x}}{2\sqrt{%
x}}
\]
This way one gets the first three Euler-Lagrange equations:
\begin{eqnarray}
\dot{r} &=&c\sqrt{x}\cos u  \nonumber \\
\dot{x} &=&-c\frac{E_{r}}{E_{x}}\sqrt{x}\cos u  \label{Lagr.reduced} \\
\dot{u} &=&-\frac{c\sin u}{2r\sqrt{x}L_{x}}\left( 2xL_{x}+rL_{r}\right)
\,\,\,\,\,\,,  \nonumber
\end{eqnarray}
which form a closed subsystem of the whole system.

Denote by $\widetilde{Z}_{L}$ the projection of $Z_{L}$ onto the $\left(
x,r,u\right) $-space. Then, solutions of (\ref{Lagr.reduced}) are identified
with trajectories of $\widetilde{Z}_{L}$. Due to obvious symmetry with
respect to the $\left( x,r\right) $-plane it is sufficient to consider those
of them for which $u\in [0,\pi ]$ (i.e. counterclockwise motions around the
center of the elastic force).

The fourth Euler-Lagrange equation
\begin{equation}
\dot{\phi}=\frac{c\sqrt{x}}{r}\sin u  \label{Lagr. quarta}
\end{equation}
can be found directly from :
\[
\tan \phi =\frac{q_{2}}{q_{1}}\,\,\,\,\,\,,
\]
so that
\begin{equation}
\phi \left( t\right) =c\int_{0}^{t}\frac{\sqrt{x\left( \tau \right) }}{%
r\left( \tau \right) }\sin u\left( \tau \right) \,d\tau
\,\,\,+const.\,\,\,\,\,\,\,\,\,\,\,,  \label{phi(t)=...}
\end{equation}
with $\left( x\left( t\right) ,r\left( t\right) ,u\left( t\right) \right) $
being a solution of (\ref{Lagr.reduced}).

Let
\begin{equation}
\Xi _{\lambda ,\mu }=\{E=\lambda mc^{2}\}\cap \{I=\mu mcr_{0}\}
\label{sigmal,m}
\end{equation}
Obviously,
\[
\Xi _{\lambda ,\mu }=\Gamma _{\lambda ,\mu }\times S^{1}\;\;,
\]
where $\Gamma _{\lambda ,\mu }$ is the projection of $\Xi _{\lambda ,\mu }$
into $\left( r,x,u\right) $-space, while the circle $S^{1}$ corresponds to
the cyclic coordinate $\phi $. In their turn surfaces $\Xi _{\lambda ,\mu }$
foliate the energy level 3-fold
\[
\Sigma _{\lambda }=\{E=\lambda mc^{2}\}
\]
In the case of oscillator (\ref{oscillator}) $\Sigma _{\lambda }$ is not
empty for $\lambda \in [1,+\infty ]$.

Let $\widetilde{\Sigma }_\lambda $ be the projection of $\Sigma _\lambda $
onto the $(r,x,u)$-space. Obviously
\[
\widetilde{\Sigma }_\lambda =\Gamma _\lambda \times S^1\,\,\,\,\,,
\]
where $\Gamma _\lambda $ is the curve in $\left( r,x\right) $-plane given by
equation $E\left( r,x\right) =\lambda mc^2$ and $S^1$ is the circle
corresponding to the cyclic coordinate $u$. Curves $\Gamma _\lambda $ are
shown in Fig \ref{fig7}. One can see that $\Gamma _\lambda $ intersects the
projections $C_1,C_2$, of $S_1,S_2$, respectively, as follows: i) at two
different points $P_1,P_2$, if $\lambda >2$; ii) at the single point $Q$, if
$\lambda =2$; iii) nowhere, if $1\leq \lambda <2$. Therefore $\Sigma
_\lambda $ intersects $S_1$ and $S_2$:

\begin{figure}
\mbox{
\epsfxsize=2.776665in
\epsfbox{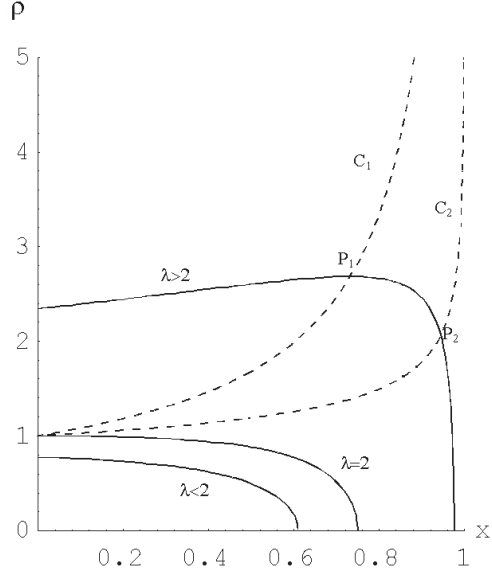}
}
\caption{Energy level surfaces}
\label{fig7}
\end{figure}

\begin{enumerate}
\item  along two tori $T_i=\Sigma _\lambda \cap S_i\,,\,i=1,2$, if $\lambda
>2$. These tori project onto $P_i$'s and have $\left( \phi ,u\right) $ as
cyclic coordinates.

\item  along the circle $\gamma =S_1\cap S_2$, with the cyclic coordinate $%
\phi $, if $\lambda =2$.

\item  nowhere if $\lambda <2$.
\end{enumerate}

\noindent Therefore, $\Sigma _{\lambda }\backslash S$ has three connected
components, if $\lambda >2$ and is connected, if $\lambda \leq 2$.

\begin{figure}
\mbox{
\epsfxsize=3.52in
\epsfbox{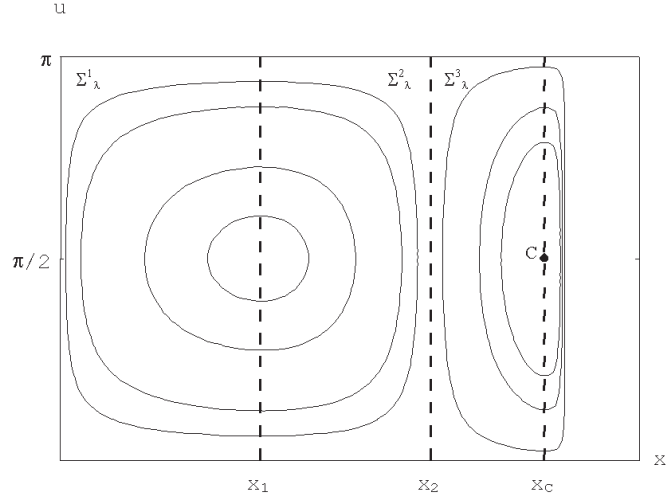}
}
\caption{Phase portrait of the oscillator for a fixed value of the energy}
\label{fig8}
\end{figure}

\begin{figure}%[p]
\mbox{
\epsfxsize=2.79in
\epsfbox{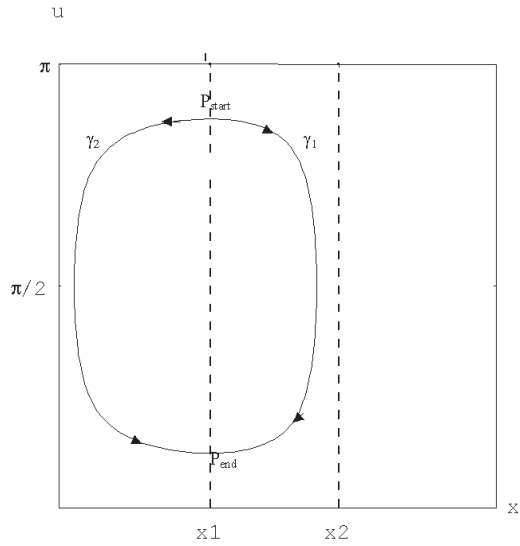}
}
\caption{Jumping phase trajectories}
\label{fig9}
\end{figure}

In the case $\lambda >2$ the behaviour of phase trajectories depends
strongly on the connected component of $\Sigma _{\lambda }\backslash S$ they
belong to. Due to (\ref{phi(t)=...}) it is sufficient to study trajectories
of (\ref{Lagr.reduced}), i.e. connected components of curves $\Gamma
_{\lambda ,\mu }\backslash \widetilde{S}$, with $\widetilde{S}$ being the
projection of $S$ onto the $\left( r,x,u\right) $-space. In Fig. \ref{fig8} the
projections into the $\left( x,u\right) $-plane of three different kinds of
such trajectories contained in $\widetilde{\Sigma }_{\lambda }$ for a fixed $%
\lambda $ are shown. As before we limit ourselves to $0\leq u\leq \pi $. The
three vertical lines correspond to the projections $\widetilde{T}_{i}$ of
tori $T_{i}$'s (in fact they are circles, due to the ciclicity of coordinate
$u$). Passing to further details, denote by $x_{i}=x_{i}\left( \lambda
\right) $ the constant value of $x$ along $T_{i}$. Connected components $%
\Sigma _{\lambda }^{1},\Sigma _{\lambda }^{2},\Sigma _{\lambda }^{3}$ of $%
\Sigma _{\lambda }\backslash S$ correspond to $x<x_{1},x_{1}<x<x_{2}$ and $%
x_{2}<x$ respectively. The trajectories belonging to $\Sigma _{\lambda }^{1}$
and $\Sigma _{\lambda }^{2}$ are discontinuous in the sense that they end on
$S_{1}$ and then jump, according to the transition principle, at another
point of $S_{1}$. More exactly, in the situation shown in the figure, such a
trajectory starts from a point of $T_{1}^{+}=T_{1}\cap \{\pi /2\leq u\leq
\pi \}$ and reaches a point in $T_{1}^{-}=T_{1}\cap \{0\leq u\leq \pi /2\}$
in a finite time. Trajectories $\gamma _{1}$ and $\gamma _{2}$, whose
projections $\widetilde{\gamma }_{1}$ and $\widetilde{\gamma }_{2}$ are
shown in Fig. \ref{fig9}, correspond to the same value of $\lambda $ and of $\mu $.
Denote by $P_{start}\equiv \left( \overline{x},\overline{r},\pi -\overline{u}%
\right) \in \widetilde{T}_{1}^{+}$ and $P_{end}\equiv (\overline{x},%
\overline{r},\overline{u})\in \widetilde{T}_{1}^{-}$ the common starting and
ending points of $\widetilde{\gamma }_{1}$ and $\widetilde{\gamma }_{2}$.
When a phase point starts from $\left( P_{start},\phi _{0}\right) \equiv
\left( \overline{x},\overline{r},\pi -\overline{u},\phi _{0}\right) $ and
then goes along $\gamma _{1}$ (or, alternatively, $\gamma _{2}$) it arrives
at the point $\left( P_{end},\overline{\phi }\right) \equiv \left( \overline{%
x},\overline{r},\overline{u},\overline{\phi }\right) $ with $\overline{\phi }
$ given by (\ref{deltaphi}) and then proceeds along the trajectory whose
projection is $\widetilde{\gamma }_{1}$ or, alternatively, $\widetilde{%
\gamma }_{2}$, and so on.

\begin{figure}%[p]
\mbox{
\epsfxsize=3.34in
\epsfbox{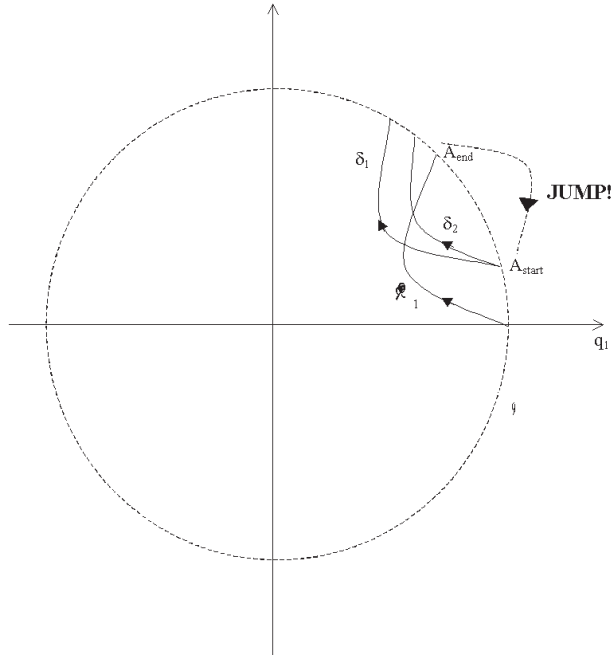}
}
\caption{Jumping oscillator}
\label{fig10}
\end{figure}

In Fig. \ref{fig10} the situation in the configuration space $M=\{\left(
q_{1},q_{2}\right) \}=\{(r,\phi )\}$ is shown.The oscillating particle,
starting from the point \linebreak
$\left( r=\overline{r},\phi =0\right) $, moves along
the projection $\widehat{\gamma }_{1}$ of $\gamma _{1}$ onto $M$ until it
reaches the point $A_{end}\equiv \left( r=\overline{r},\phi =\overline{\phi }%
\right) $. From there it jumps to the point $A_{start}\acute{\equiv}\left( r=%
\overline{r},\phi =\widetilde{\phi }\right) $, where $\widetilde{\phi }$ is
given by (\ref{deltaphi}). Then it splits into the two trajectories $\delta
_{1},\delta _{2}$, both with initial velocity
$$
\left( x=\overline{x},u=%
\widetilde{u}=\pi -\overline{u}\right).
$$

\begin{figure}
\mbox{
\epsfxsize=2.11in
\epsfbox{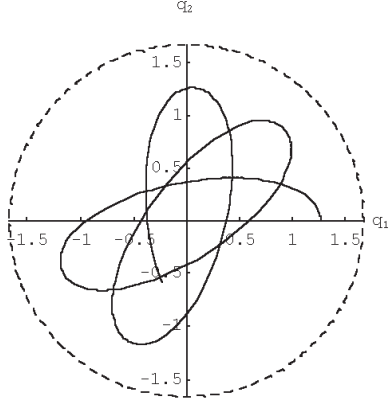}
}
\caption{Regular precessions}
\label{fig11}
\end{figure}

Consider now the trajectories of (\ref{Lagr.reduced}) contained in $\Sigma
_{\lambda }^{3}$. These are regular closed trajectories winding around the
center $C=C\left( \lambda \right) \equiv \left( x=x_{C},u=\pi /2\right) $,
where $x_{C}=x_{C}\left( \lambda \right) $ is the zero of the equation $%
2xL_{x}\left( r_{en}\left( x,\lambda \right) ,x\right) +rL_{r}\left(
r_{en}\left( x,\lambda \right) ,x\right) =0$, and $r=r_{en}\left( x,\lambda
\right) $ is implicitly defined by the equation $E\left( r,x\right) =\lambda
mc^{2}$. The corresponding trajectories in the configuration space are shown
in Fig. \ref{fig11}. They are precessions around the force center.

Note that both the discontinuous trajectories in $\Sigma _{\lambda }^{2}$
and the regular ones in $\Sigma _{\lambda }^{3}$ near $S_{2}$ tend to be
parallel to it, so that this component of the singular hypersurface is never
reached by the phase point, at least for trajectories with non-zero angular
momentum.

\medskip

\textbf{Aknowledgments}. The second author is obliged to Gaetano Vilasi for
drawing his attention to not everywhere regular actions. His thanks go also
to Igor Pavlotsky for a discussion of his recent results concerning Darwin's
model.

\end{document}